\documentclass{article}%
\usepackage{amsfonts}
\usepackage{amssymb}
\usepackage{amsmath}
\usepackage{graphicx}%
\providecommand{\U}[1]{\protect\rule{.1in}{.1in}}
\newtheorem{theorem}{Theorem}

\newtheorem{example}[theorem]{Example}

\begin{document}

\title{A Graph-theoretic Method to Define any Boolean Operation on Partitions}
\author{David Ellerman\\Philosophy Department,\\University of California, \\Riverside, CA 92521 USA\\and\\School of Social Sciences,\\University of Ljubljana, \\Slovenia\\david@ellerman.org\\orcid.org/0000-0002-5718-618X\\Published as: \textquotedblleft A Graph-Theoretic Method to Define \\Any Boolean Operation on Partitions.\textquotedblright\ \\\textit{The Art of Discrete and Applied Mathematics} 2 (2): 1--9.\\https://doi.org/10.26493/2590-9770.1259.9d5.}
\maketitle

\begin{abstract}
\noindent The lattice operations of join and meet were defined for set
partitions in the nineteenth century, but no new logical operations on
partitions were defined and studied during the twentieth century. Yet there is
a simple and natural graph-theoretic method presented here to define any
$n$-ary Boolean operation on partitions. An equivalent closure-theoretic
method is also defined. In closing, the question is addressed of why it took
so long for all Boolean operations to be defined for partitions.



\end{abstract}

\section{Introduction}

The lattice operations of join and meet were defined on set partitions during
the late nineteenth century, and the lattice of partitions on a set was used
as an example of a non-distributive lattice. But during the entire twentieth
century, no new logical operations were defined on partitions.

\begin{quotation}
Equivalence relations are so ubiquitous in everyday life that we often forget
about their proactive existence. Much is still unknown about equivalence
relations. Were this situation remedied, the theory of equivalence relations
could initiate a chain reaction generating new insights and discoveries in
many fields dependent upon it.

This paper springs from a simple acknowledgement: the only operations on the
family of equivalence relations fully studied, understood and deployed are the
binary join $\vee$ and meet $\wedge$ operations. \cite[p. 445]{bmp:eqrel}
\end{quotation}

\noindent Papers on the "logic" of equivalence relations \cite{rota:logiccers}
or partitions only involved the join and meet, and not the crucial logical
operation of implication.

Yet, there is a general graph-theoretic method\footnote{The method is,
strictly speaking, an algorithm only when $U$ is finite.} by which any $n$-ary
Boolean (or truth-functional) operation $f:\left\{  T,F\right\}
^{n}\rightarrow\left\{  T,F\right\}  $ can be used to define the corresponding
$n$-ary operation $f:\prod\left(  U\right)  ^{n}\rightarrow\prod\left(
U\right)  $ where $\prod\left(  U\right)  $ is the set of partitions on a set
$U$.

A \textit{partition} $\pi=\left\{  B,B^{\prime},...\right\}  $ on a set
$U=\left\{  u,u^{\prime},...\right\}  $ is a set of disjoint non-empty subsets
$B,B^{\prime},...$ of $U$, called \textit{blocks}, whose union is $U$. The
corresponding equivalence relation, denoted $\operatorname*{indit}\left(
\pi\right)  $, is the set of ordered pairs of elements of $U$ that are in the
same block of $\pi$, and are called the \textit{indistinctions}
or\textit{\ indits} of $\pi$, i.e.,

\begin{center}
$\operatorname*{indit}\left(  \pi\right)  =\left\{  \left(  u,u^{\prime
}\right)  \in U\times U:\exists B\in\pi,u,u^{\prime}\in B\right\}  $.
\end{center}

\noindent The complement $\operatorname*{dit}\left(  \pi\right)  =U\times
U-\operatorname*{indit}\left(  \pi\right)  $ is the set of
\textit{distinctions} or \textit{dits} of $\pi$, i.e., ordered pairs of
elements in different blocks. As binary relations, the sets of distinctions or
\textit{ditsets} $\operatorname*{dit}\left(  \pi\right)  $ of some partition
$\pi$ on $U$ are called\textit{\ partition} (or \textit{apartness})
\textit{relations.} Given partitions $\pi=\left\{  B,B^{\prime},...\right\}  $
and $\sigma=\left\{  C,C^{\prime},...\right\}  $ on $U$, the
\textit{refinement} relation is the partial order defined by:

\begin{center}
$\sigma\preceq\pi$ if $\forall B\in\pi,\exists C\in\sigma,B\subseteq C $.
\end{center}

\noindent At the top of the refinement partial order is the \textit{discrete
partition} $\mathbf{1}=\left\{  \left\{  u\right\}  :u\in U\right\}  $ of all
singletons and at the bottom is the \textit{indiscrete partition}
$\mathbf{0}=\left\{  U\right\}  $ with only one block consisting of $U$. In
terms of binary relations, the refinement partial order is just the inclusion
partial order on ditsets, i.e., $\sigma\preceq\pi$ iff $\operatorname*{dit}%
\left(  \sigma\right)  \subseteq\operatorname*{dit}\left(  \pi\right)  $. It
should be noted that most of the previous literature on partitions (e.g.,
\cite{Birk:lt}) uses the opposite partial order of `unrefinement'
corresponding to the inclusion relation on equivalence relations--which
reverses the definitions of the join and meet of partitions.

\section{The Join Operation on Partitions}

The \textit{join} $\pi\vee\sigma$ of partitions $\pi$ and $\sigma$ (least
upper bound using the refinement partial order) is the partition whose blocks
are the non-empty intersections $B\cap C$ of the blocks of $\pi$ and $\sigma$
(under the unrefinement ordering, it is the meet). In terms of ditsets,
$\operatorname*{dit}\left(  \pi\vee\sigma\right)  =\operatorname*{dit}\left(
\pi\right)  \cup\operatorname*{dit}\left(  \sigma\right)  $. The general
method for defining Boolean operations on partitions will be first illustrated
with the join operation whose corresponding Boolean operation is disjunction
with the truth table.

\begin{center}%
\begin{tabular}
[c]{|c|c|c|}\hline
$P$ & $Q$ & $P\vee Q$\\\hline\hline
$T$ & $T$ & $T$\\\hline
$T$ & $F$ & $T$\\\hline
$F$ & $T$ & $T$\\\hline
$F$ & $F$ & $F$\\\hline
\end{tabular}

Truth table for disjunction.
\end{center}

Let $K\left(  U\right)  $ be the complete undirected graph on $U$. The links
$u-u^{\prime}$ corresponding to dits, i.e., $\left(  u,u^{\prime}\right)
\in\operatorname*{dit}\left(  \pi\right)  $, of a partition are labelled with
the `truth value' $T_{\pi}$ and corresponding to indits $\left(  u,u^{\prime
}\right)  \in\operatorname*{indit}\left(  \pi\right)  $ are labelled with the
`truth value' $F_{\pi}$. Given the two partitions $\pi$ and $\sigma$, each
link in the complete graph $K\left(  U\right)  $ is labelled with a pair of
truth values. The \textit{graph }$G\left(  \pi\vee\sigma\right)  $\textit{\ of
the join }is obtained by putting a link $u-u^{\prime}$ where the truth
function applied to the pair of truth values on the link in $K\left(
U\right)  $ gives an $F$. Thus in the case at hand, the only links in
$G\left(  \pi\vee\sigma\right)  $ are for the $u-u^{\prime}$ labelled with
$F_{\pi}$ and $F_{\sigma}$ in $K\left(  U\right)  $. Then the partition
$\pi\vee\sigma$ is obtained as the connected components of its graph $G\left(
\pi\vee\sigma\right)  $. Thus $u$ and $u^{\prime}$ are in the same block
(connected component of $G\left(  \pi\vee\sigma\right)  $) if and only if the
link $u-u^{\prime}$ was labelled $F_{\pi}$ and $F_{\sigma}$, i.e., $u$ and
$u^{\prime}$ were in the same block of $\pi$ and in the same block of $\sigma
$. Thus the graph-theoretic definition of the join reproduces the
set-of-blocks definition of the join defined as having its blocks the
non-empty intersections of the blocks of $\pi$ and $\sigma$.

\section{The Meet Operation on Partitions}

On the combined set of blocks $\pi\cup\sigma$ of $\pi$ and $\sigma$, define
the \textit{overlap relation} $B\between C$ on two blocks if they have a
non-empty intersection or overlap (see \cite{ore:ter}). The
reflexive-symmetric-transitive closure of this relation is an equivalence
relation, and the union of the blocks in each equivalence class gives the
blocks of the \textit{meet} $\pi\wedge\sigma$. The corresponding
truth-functional operation is conjunction with the following truth table.

\begin{center}%
\begin{tabular}
[c]{|c|c|c|}\hline
$P$ & $Q$ & $P\wedge Q$\\\hline\hline
$T$ & $T$ & $T$\\\hline
$T$ & $F$ & $F$\\\hline
$F$ & $T$ & $F$\\\hline
$F$ & $F$ & $F$\\\hline
\end{tabular}

Truth table for conjunction.
\end{center}

The same method is applied except that the links of the graph $G\left(
\pi\wedge\sigma\right)  $ are the ones for which the conjunction truth table
gives an $F$ when applied to the truth values on each link $u-u^{\prime}$.
Thus $G\left(  \pi\wedge\sigma\right)  $ contains a link $u-u^{\prime}$ if
$\left(  u,u^{\prime}\right)  \in\operatorname*{indit}\left(  \pi\right)  $,
$\left(  u,u^{\prime}\right)  \in\operatorname*{indit}\left(  \sigma\right)
$, or both. Then the blocks of the partition $\pi\wedge\sigma$ are the
connected components of the graph $G\left(  \pi\wedge\sigma\right)  $.

The proof that the graph-theoretic definition of the meet gives the usual
set-of-blocks definition of the meet boils down to showing that: $B\in\pi$ and
$C\in\sigma$ are contained in the same block of the usual meet $\pi
\wedge\sigma$ (i.e., there is a chain of overlaps $B\between C^{\prime
}\between...\between B^{\prime}\between C$ connecting $B$ and $C$) if and only
for any $u\in B$ and $u^{\prime}\in C$, $u$ and $u^{\prime}$ are in the same
connected component of $G\left(  \pi\wedge\sigma\right)  $. If any two blocks
$B^{\prime}\between C^{\prime}$ overlap in the overlap chain, then there is an
element $u^{\prime\prime}\in B^{\prime}\cap C^{\prime}$ such any $u\in
B^{\prime}$ had a link $u-u^{\prime\prime}$ in $G\left(  \pi\wedge
\sigma\right)  $ and similarly any $u^{\prime}\in C^{\prime}$ has a link
$u^{\prime\prime}-u^{\prime}$ in $G\left(  \pi\wedge\sigma\right)  $. Hence
the existence of an overlap chain connecting $B$ and $C$ implies that any
$u\in B$ and $u^{\prime}\in C$ are in the same connected component of
$G\left(  \pi\wedge\sigma\right)  $. Conversely, if $u\in B$ and $u^{\prime
}\in C$ are in the same connected component of $G\left(  \pi\wedge
\sigma\right)  $, then there is some chain of links $u=u_{0}-u_{1}%
-...-u_{n-1}-u_{n}=u^{\prime}$ where each link $u_{i}-u_{i+1}$ for
$i=0,...,n-1$ has either $\left(  u_{i},u_{i+1}\right)  \in
\operatorname*{indit}\left(  \pi\right)  $, $\left(  u_{i},u_{i+1}\right)
\in\operatorname*{indit}\left(  \sigma\right)  $, or both. Every link
$u_{i}-u_{i+1}$ that is in one indit set but not the other, say, $\left(
u_{i},u_{i+1}\right)  \in\operatorname*{indit}\left(  \pi\right)  $ and
$\left(  u_{i},u_{i+1}\right)  \notin\operatorname*{indit}\left(
\sigma\right)  $, establishes an overlap between the block of $\pi$ containing
$u_{i},u_{i+1}$ and the block of $\sigma$ containing $u_{i}$ as well as the
different block of $\sigma$ containing $u_{i+1}$. Thus the chain of links
connecting $u\in B$ and $u^{\prime}\in C$ establishes a chain of overlapping
blocks connecting $B$ and $C$.

\section{The Implication Operation on Partitions}

The real beginning of the \textit{logic} of partitions, as opposed to the
lattice theory of partitions, was the discovery of the set-of-blocks
definition of the implication operation $\sigma\Rightarrow\pi$ for partitions
(\cite{ell:lop}, \cite{ell:intropartlogic}). The intuitive idea is that
$\sigma\Rightarrow\pi$ functions like an indicator or characteristic function
to indicate which blocks $B$ of $\pi$ are contained in a block of $\sigma$.
View the discretized version of $B\in\pi$, i.e., $B$ replaced by the set of
singletons of the elements of $B$, as the local version $\mathbf{1}_{B}$ of
the discrete partition $\mathbf{1}$, and view the block $B$ remaining whole as
the local version $\mathbf{0}_{B}$ of the indiscrete partition $\mathbf{0}$.
Then the partition implication as the inclusion indicator function is: the
blocks of $\sigma\Rightarrow\pi$ are for any $B\in\pi$:

\begin{center}
$\left\{
\begin{array}
[c]{l}%
\mathbf{1}_{B}\text{ if }\exists C\in\sigma,B\subseteq C\\
\mathbf{0}_{B}=B\text{ otherwise.}%
\end{array}
\right.  .$
\end{center}

In the case of the Boolean logic of subsets, for any subsets $S,T\subseteq U$,
the conditional $S\supset T=S^{c}\cup T$ has the property: $S\supset T=U$ iff
$S\subseteq T$, i.e., the conditional $S\supset T$ equals the top of the
lattice of subsets of $U$ iff the inclusion relation $S\subseteq T$ holds.
Similarly, it is immediate that the corresponding relation holds in the
partition case:

\begin{center}
$\sigma\Rightarrow\pi=\mathbf{1}$ iff $\sigma\preceq\pi$.
\end{center}

\noindent This set-of-blocks definition of the partition implication operation
accounts for the important new non-lattice-theoretic properties revealed in
the \textit{algebra} of partitions $\prod\left(  U\right)  $ on $U $ (defined
with the join, meet, and implication as partition operations).

A logical formula in the language of join, meet, and implication is a
\textit{subset tautology} if for any non-empty universe $U$ and any subsets of
$U$ substituted for the variables, the whole formula evaluates by the
set-theoretic operations of join, meet, and implication (conditional) to the
top $U$. Similarly, a formula in the same language is a \textit{partition
tautology} if for any universe $U$ with $\left\vert U\right\vert >1$ and for
any partitions on $U$ substituted for the variables, the whole formula
evaluates by the partition operations of join, meet, and implication to the
top $\mathbf{1}$ (the discrete partition). All partition tautologies are
subset tautologies but not vice-versa. \textit{Modus ponens} $\left(
\sigma\wedge\left(  \sigma\Rightarrow\pi\right)  \right)  \Rightarrow\pi$ is
both a subset and partition tautology but Peirce's law, $\left(  \left(
\sigma\Rightarrow\pi\right)  \Rightarrow\sigma\right)  \Rightarrow\sigma$,
accumulation, $\sigma\Rightarrow\left(  \pi\Rightarrow\left(  \sigma\wedge
\pi\right)  \right)  $, and distributivity, $\left(  \left(  \pi\vee
\sigma\right)  \wedge\left(  \pi\vee\tau\right)  \right)  \Rightarrow\left(
\pi\vee\left(  \sigma\wedge\tau\right)  \right)  $, are examples of subset
tautologies that are not partition tautologies. The importance of the
implication for partition logic is emphasized by the fact that the only
partition tautologies using only the lattice operations, e.g., $\pi
\vee\mathbf{1}$, correspond to general lattice-theoretic identities, i.e.,
$\pi\vee\mathbf{1}=\mathbf{1}$ (see \cite{sachs:iden}).

The graph-theoretic method automatically gives a partition operation
corresponding to the Boolean conditional or implication with the truth table:

\begin{center}%
\begin{tabular}
[c]{|c|c|c|}\hline
$P$ & $Q$ & $P\supset Q$\\\hline\hline
$T$ & $T$ & $T$\\\hline
$T$ & $F$ & $F$\\\hline
$F$ & $T$ & $T$\\\hline
$F$ & $F$ & $T$\\\hline
\end{tabular}

Truth table for conditional
\end{center}

\noindent and it is not trivial that the two definitions are the same. It may
be helpful to restate the truth table in terms of the partitions.

\begin{center}%
\begin{tabular}
[c]{|c|c|c|}\hline
$\sigma$ & $\pi$ & $\sigma\Rightarrow\pi$\\\hline\hline
$T_{\sigma}$ & $T_{\pi}$ & $T_{\sigma\Rightarrow\pi}$\\\hline
$T_{\sigma}$ & $F_{\pi}$ & $F_{\sigma\Rightarrow\pi}$\\\hline
$F_{\sigma}$ & $T_{\pi}$ & $T_{\sigma\Rightarrow\pi}$\\\hline
$F_{\sigma}$ & $F_{\pi}$ & $T_{\sigma\Rightarrow\pi}$\\\hline
\end{tabular}

Implication truth table for partition `truth values'.
\end{center}

For the graph-theoretic definition of $\sigma\Rightarrow\pi$, we again label
the links $u-u^{\prime}$ in the complete graph $K\left(  U\right)  $ with
$T_{\pi}$ if $\left(  u,u^{\prime}\right)  \in\operatorname*{dit}\left(
\pi\right)  $ and $F_{\pi}$ otherwise, and similarly for $\sigma$. Then we
construct the graph $G\left(  \sigma\Rightarrow\pi\right)  $ by putting in a
link $u-u^{\prime}$ only in the case the link is labeled $T_{\sigma}$ and
$F_{\pi}$, i.e., $F_{\sigma\Rightarrow\pi}$. Then the partition $\sigma
\Rightarrow\pi$ is the partition of connected components in the graph
$G\left(  \sigma\Rightarrow\pi\right)  $.

To prove the graph-theoretic and set-of-blocks definitions equivalent, we
might first note that if $\left(  u,u^{\prime}\right)  \in\operatorname*{dit}%
\left(  \pi\right)  $, then $T_{\pi}$ is assigned to that link in $K\left(
U\right)  $ so there is no link $u-u^{\prime}$ in $G\left(  \sigma
\Rightarrow\pi\right)  $. And if $\left(  u,u^{\prime}\right)  \in
\operatorname*{indit}\left(  \pi\right)  $ but also $\left(  u,u^{\prime
}\right)  \in\operatorname*{indit}\left(  \sigma\right)  $, then
$T_{\sigma\Rightarrow\pi}$ is assigned to the link in $K\left(  U\right)  $ so
again there is no link $u-u^{\prime}$ in $G\left(  \sigma\Rightarrow
\pi\right)  $. There is a link $u-u^{\prime}$ in $G\left(  \sigma
\Rightarrow\pi\right)  $ in and only in the following situation where $\left(
u,u^{\prime}\right)  \in\operatorname*{indit}\left(  \pi\right)  $ and
$\left(  u,u^{\prime}\right)  \in\operatorname*{dit}\left(  \sigma\right)
$--which is exactly the situation when $B$ is not contained in any block $C$
of $\sigma$:%

%TCIMACRO{\FRAME{dtbpF}{2.3713in}{1.081in}{0pt}{}{}{fig1-sigmatruepifalse.eps}%
%width 2.3713in;  height 1.081in;  depth 0pt;  original-width 2.3307in;
%original-height 1.0482in;  cropleft "0";  croptop "1";  cropright "1";
%filename 'fig1-sigmatruepifalse.eps';file-properties "XNPEU";}} }%
%BeginExpansion
\begin{center}
\includegraphics[
height=1.081in,
width=2.3713in
]%
{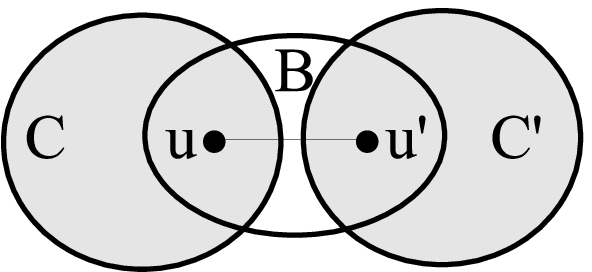}%
\end{center}
%EndExpansion

\begin{center}
Figure 1: Links $u-u^{\prime}$ in $G\left(  \sigma\Rightarrow\pi\right)  $.
\end{center}

\noindent Then for any other element $u^{\prime\prime}\in B$ so that $\left(
u,u^{\prime\prime}\right)  $ and $\left(  u^{\prime},u^{\prime\prime}\right)
\in\operatorname*{indit}\left(  \pi\right)  $, we must have either $\left(
u,u^{\prime\prime}\right)  \in\operatorname*{dit}\left(  \sigma\right)  $ or
$\left(  u^{\prime},u^{\prime\prime}\right)  \in\operatorname*{dit}\left(
\sigma\right)  $ so $u^{\prime\prime}$ is linked in $G\left(  \sigma
\Rightarrow\pi\right)  $ to either $u$ or to $u^{\prime}$. Thus all the
elements of $B$ are in the same connected component of the graph $G\left(
\sigma\Rightarrow\pi\right)  $ whenever $B$ is not contained in any block of
$\sigma$. If, on the other hand, $B$ is contained in some block $C$ of
$\sigma$, then any $u\in B$ cannot be linked to any other $u^{\prime}$. In
order to that $F_{\pi}$ assigned to the link $u-u^{\prime}$, the two elements
have to both belong to $B$ and thus since $B\subseteq C$, they both belong to
$C$ so $F_{\sigma}$ and thus $T_{\sigma\Rightarrow\pi}$ is also assigned to
that link. Thus when $B$ is contained in a block $C\in\sigma$, then any point
$u\in B$ is a disconnected component to itself in $G\left(  \sigma
\Rightarrow\pi\right)  $ so $B$ is discretized in the graph-theoretic
construction of $\sigma\Rightarrow\pi$. Thus the graph-theoretic and
set-of-blocks definitions of the partition implication are equivalent.%

%TCIMACRO{\FRAME{dtbpF}{4.2566in}{1.8723in}{0pt}{}{}%
%{fig2-graph4implication.eps}{\special{ language "Scientific Word";
%type "GRAPHIC";  maintain-aspect-ratio TRUE;  display "USEDEF";
%valid_file "F";  width 4.2566in;  height 1.8723in;  depth 0pt;
%original-width 5.3454in;  original-height 2.3359in;  cropleft "0";
%filename 'fig2-graph4implication.eps';file-properties "XNPEU";}} }%
%BeginExpansion
\begin{center}
\includegraphics[
height=1.8723in,
width=4.2566in
]%
{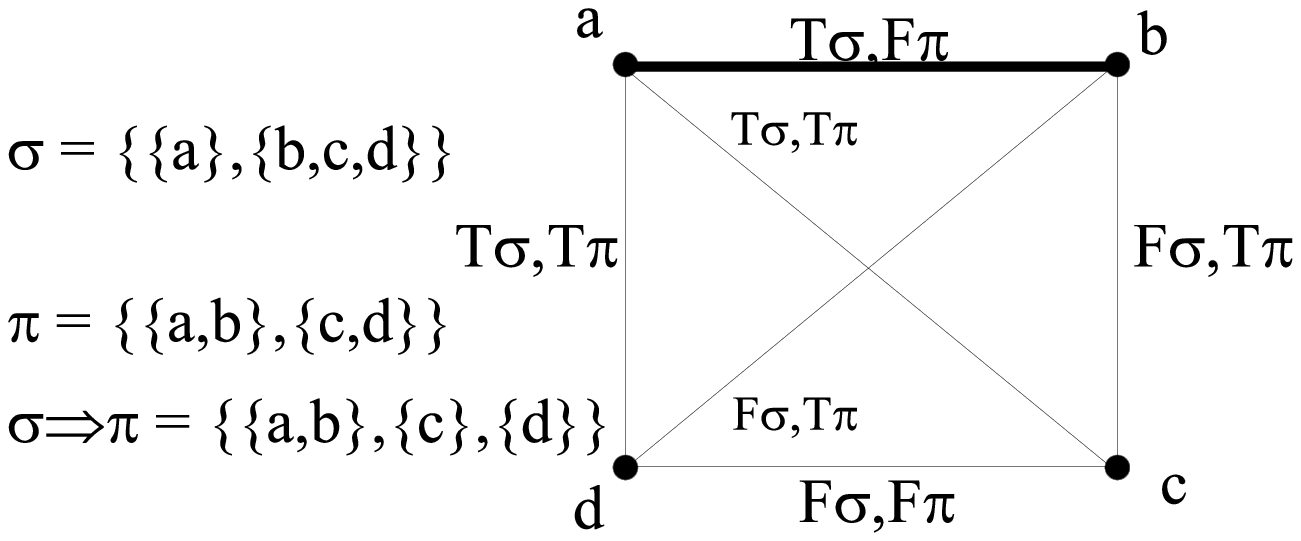}%
\end{center}
%EndExpansion

\begin{center}
Figure 2: Example of graph for partition implication
\end{center}

\begin{example}
Let $U=\left\{  a,b,c,d\right\}  $ so that $K(U)=K_{4}$ is the complete graph
on four points. Let $\sigma=\left\{  \left\{  a\right\}  ,\left\{
b,c,d\right\}  \right\}  $ and $\pi=\left\{  \left\{  a,b\right\}  ,\left\{
c,d\right\}  \right\}  $ so we see immediately from the set-of-blocks
definition, that the $\pi$-block of $\left\{  c,d\right\}  $ will be
discretized while the $\pi$-block of $\left\{  a,b\right\}  $ will remain
whole so the partition implication is $\sigma\Rightarrow\pi=\left\{  \left\{
a,b\right\}  ,\left\{  c\right\}  ,\left\{  d\right\}  \right\}  $. After
labelling the links in $K\left(  U\right)  $, we see that only the $a-b$ link
has the $F_{\sigma\Rightarrow\pi}$ `truth value' so the graph $G\left(
\sigma\Rightarrow\pi\right)  $ has only that $a-b$ link (thickened in Figure
2). Then the connected components of $G\left(  \sigma\Rightarrow\pi\right)  $
give the same partition implication $\sigma\Rightarrow\pi=\left\{  \left\{
a,b\right\}  ,\left\{  c\right\}  ,\left\{  d\right\}  \right\}  $.
\end{example}

The partition implication is quite rich in defining new structures in the
algebra of partitions (i.e., the lattice of partitions extended with other
partition operations such as the implication). For instance, for a fixed
partition $\pi$ on $U$, all the partitions of the form $\sigma\Rightarrow\pi$
(for any partitions $\sigma$ on $U$) form a Boolean algebra under the
partition operations of implication, join, and meet, e.g., $\left(
\sigma\Rightarrow\pi\right)  \Rightarrow\pi$ is the negation of $\sigma
\Rightarrow\pi$, called the \textit{Boolean core} of the upper segment
$\left[  \pi,\mathbf{1}\right]  $ in the partition algebra $\prod\left(
U\right)  $.

A \textit{relation} is a subset of a product, and, dually, a
\textit{corelation} is a partition on a coproduct. Any partition $\pi$ on $U$
can be canonically represented as a relation: $\operatorname*{dit}\left(
\pi\right)  \subseteq U\times U$. Dually any subset $S\subseteq U$ can be
canonically represented as a corelation, namely the partition $\pi\left(
S\right)  $ on the coproduct (disjoint union) $U\uplus U$ where the only
nonsingleton blocks in $\pi\left(  S\right)  $ are the pairs $\left\{
u,u^{\ast}\right\}  $ of $u$ and its copy $u^{\ast}$ for $u\notin S$. Using
this corelation construction, any powerset Boolean algebra $\wp\left(
U\right)  $ can be canonically represented as the Boolean core of the upper
segment $\left[  \pi,\mathbf{1}\right]  $ in the partition algebra
$\prod\left(  U\uplus U\right)  $ where $\pi=\pi\left(  \emptyset\right)  $ is
the partition on the disjoint union $U\uplus U$ whose blocks are all the pairs
$\left\{  u,u^{\ast}\right\}  $ for each element $u\in U$ and its copy
$u^{\ast}$. Each partition of the form $\sigma\Rightarrow\pi$ on $U\uplus U$
is $\pi\left(  S\right)  $ for some $S\subseteq U$ since $\sigma\Rightarrow
\pi$ is essentially the characteristic function of some subset $S$ of $U$ with
$\mathbf{1}\Rightarrow\pi=\pi\left(  \emptyset\right)  $ playing the role of
the empty set $\emptyset$ and $\pi\Rightarrow\pi=\mathbf{1}_{U\uplus U}$
playing the role of $U$.

\section{The General Graph-Theoretic Method}

Let $f:\left\{  T,F\right\}  ^{n}\rightarrow\left\{  T,F\right\}  $ be an
$n$-ary Boolean function and let $\pi_{1},...,\pi_{n}$ be $n$ partitions on
$U$. In order to define the corresponding $n$-ary partition operation
$f\left(  \pi_{1},...,\pi_{n}\right)  $, we again consider the complete graph
$K\left(  U\right)  $ and then use each partition $\pi_{i}$ to label each link
$u-u^{\prime}$ with $T_{\pi_{i}}$ if $\left(  u,u^{\prime}\right)
\in\operatorname*{dit}\left(  \pi_{i}\right)  $ and $F_{\pi_{i}}$ if $\left(
u,u^{\prime}\right)  \in\operatorname*{indit}\left(  \pi_{i}\right)  $. Then
on each link we may apply $f$ to the $n$ `truth values' on the link and retain
the link in $G\left(  f\left(  \pi_{1},...,\pi_{n}\right)  \right)  $ if the
result was $F_{f\left(  \pi_{1},...,\pi_{n}\right)  }$. The partition
$f\left(  \pi_{1},...,\pi_{n}\right)  $ is obtained as the connected
components of the graph $G\left(  f\left(  \pi_{1},...,\pi_{n}\right)
\right)  $.

\section{An Equivalent Closure-theoretic Method}

Given any subset $S\subseteq U\times U$,
the\textit{\ reflexive-symmetric-transitive (RST) closure} $\overline{S}$ is
the intersection of all equivalence relations on $U$ containing $S$. The
`topological' terminology of calling a subset \textit{closed} if
$S=\overline{S}$ is used even though the RST closure operator is \textit{not}
a topological closure operator since the union of two closed sets is not
necessarily closed. The closed sets in $U\times U$ are the equivalence
relations (or indit sets of partitions), and their complements, the
\textit{open} sets, are the partition relations (or ditsets of partitions). As
usual, the interior operator $\operatorname*{int}\left(  S\right)  =\left(
\overline{S^{c}}\right)  ^{c}$ is the complement of the closure of the
complement, and the open sets are the ones equalling their interiors.

The closure-theoretic method of defining Boolean operations on partitions will
be illustrated using the symmetric difference or inequivalence operation
$\pi\oplus\sigma$. Every $n$-ary Boolean operation can be defined by a truth
table such as the one for symmetric difference in this case:

\begin{center}%
\begin{tabular}
[c]{|c|c|c|}\hline
$P$ & $Q$ & $P\oplus Q$\\\hline\hline
$T$ & $T$ & $F$\\\hline
$T$ & $F$ & $T$\\\hline
$F$ & $T$ & $T$\\\hline
$F$ & $F$ & $F$\\\hline
\end{tabular}

Truth table for symmetric difference.
\end{center}

The disjunctive normal form (DNF) for the formula $P\oplus Q$ is given by the
rows where the formula evaluates as $T$, i.e., $P\oplus Q=\left(  P\wedge\lnot
Q\right)  \vee\left(  \lnot P\wedge Q\right)  $, while the DNF for the
negation of the formula is given by the other rows where the formula evaluates
as $F$, i.e., $\lnot\left(  P\oplus Q\right)  =\left(  P\wedge Q\right)
\vee\left(  \lnot P\wedge\lnot Q\right)  $. Given two partitions $\pi$ and
$\sigma$ on $U$, the closure-theoretic method of obtaining the partition
$\pi\oplus\sigma$ is to start with the DNF for the negated Boolean formula and
replace each unnegated variable by the corresponding ditset and each negated
variable by the corresponding indit set--as well as replacing the disjunctions
and conjunctions by the corresponding subset operations of union and
intersection. Applied to $\lnot\left(  P\oplus Q\right)  =\left(  P\wedge
Q\right)  \vee\left(  \lnot P\wedge\lnot Q\right)  $, this procedure would
yield $\left(  \operatorname*{dit}\left(  \pi\right)  \cap\operatorname*{dit}%
\left(  \sigma\right)  \right)  \cup\left(  \operatorname*{indit}\left(
\pi\right)  \cap\operatorname*{indit}\left(  \sigma\right)  \right)  \subseteq
U\times U$. Then the indit set of $\pi\oplus\sigma$ is obtained as the RST closure:

\begin{center}
$\operatorname*{indit}\left(  \pi\oplus\sigma\right)  =\overline{\left(
\operatorname*{dit}\left(  \pi\right)  \cap\operatorname*{dit}\left(
\sigma\right)  \right)  \cup\left(  \operatorname*{indit}\left(  \pi\right)
\cap\operatorname*{indit}\left(  \sigma\right)  \right)  }$
\end{center}

\noindent and the partition $\pi\oplus\sigma$ is the set of equivalence
classes of this equivalence relation.

The graph-theoretic method of obtaining the partition $\pi\oplus\sigma$ would
label each link $u-u^{\prime}$ in $K\left(  U\right)  $ by the two `truth
values' given by $\pi$ and $\sigma$, and then retain in the graph $G\left(
\pi\oplus\sigma\right)  $ the links where the truth values evaluated to
$F_{\pi\oplus\sigma}$, namely the ones labelled with $T_{\pi},T_{\sigma}$ and
$F_{\pi},F_{\sigma}$. Then the partition $\pi\oplus\sigma$ is obtained as the
connected components of the graph $G\left(  \pi\oplus\sigma\right)  $.

To see the equivalence between the two methods, note first that the links
retained in $G\left(  \pi\oplus\sigma\right)  $ are precisely the pairs
$\left(  u,u^{\prime}\right)  $ in $\left(  \operatorname*{dit}\left(
\pi\right)  \cap\operatorname*{dit}\left(  \sigma\right)  \right)  \cup\left(
\operatorname*{indit}\left(  \pi\right)  \cap\operatorname*{indit}\left(
\sigma\right)  \right)  $. The equivalence proof is completed by showing that
taking connected components in the graph $G\left(  \pi\oplus\sigma\right)  $
is equivalent to taking the RST closure of $\left(  \operatorname*{dit}\left(
\pi\right)  \cap\operatorname*{dit}\left(  \sigma\right)  \right)  \cup\left(
\operatorname*{indit}\left(  \pi\right)  \cap\operatorname*{indit}\left(
\sigma\right)  \right)  $. The elements $u$ and $u^{\prime}$ are in the same
connected component of $G\left(  \pi\oplus\sigma\right)  $ iff there is a
chain of links $u=u_{0}-u_{1}-...-u_{n-1}-u_{n}=u^{\prime}$ in the graph
$G\left(  \pi\oplus\sigma\right)  $ so each link has to be originally labelled
$T_{\pi},T_{\sigma}$ or $F_{\pi},F_{\sigma}$ in the graph on $K\left(
U\right)  $. But the condition for $\left(  u,u^{\prime}\right)  $ to be
included in the RST closure $\overline{\left(  \operatorname*{dit}\left(
\pi\right)  \cap\operatorname*{dit}\left(  \sigma\right)  \right)  \cup\left(
\operatorname*{indit}\left(  \pi\right)  \cap\operatorname*{indit}\left(
\sigma\right)  \right)  }$ is that there is a chain of pairs $\left(
u,u_{1}\right)  ,\left(  u_{1},u_{2}\right)  ,...,\left(  u_{n-1},u^{\prime
}\right)  $ such that each pair is either in $\operatorname*{dit}\left(
\pi\right)  \cap\operatorname*{dit}\left(  \sigma\right)  $ or in
$\operatorname*{indit}\left(  \pi\right)  \cap\operatorname*{indit}\left(
\sigma\right)  $. Hence the two methods give the same result.

The example suffices to illustrate the general closure-theoretic method and
its equivalence to the graph-theoretic method of defining Boolean operations
on partitions.

\section{Relationships between Boolean operations on partitions}

For two subset variables, there are $2^{4}=16$ binary Boolean operations on
subsets--corresponding to the sixteen ways to fill in the truth table for a
binary Boolean operation. Any compound Boolean function of two variables will
be truth-table equivalent to one of the sixteen binary Boolean operations. For
instance, the Pierce's Law formula $\left(  \left(  Q\Rightarrow P\right)
\Rightarrow Q\right)  \Rightarrow Q$ defines a compound binary operation that
is equivalent to the constant function $T$ since it is a subset tautology.
Certain subsets of the sixteen binary operations suffice to define all the
binary operations, e.g., $\lnot$ and $\vee$.

Matters are rather different for the Boolean operations on partitions. Using
the graph-theoretic or the closure-theoretic method, partition versions of
sixteen binary Boolean operations are easily defined. And certain combinations
of the sixteen operations suffice to define all sixteen, e.g., $\vee$,
$\wedge$, $\Rightarrow$, and $\oplus$ \cite[309-310 and fn. 18]{ell:lop}. But
when the sixteen operations are compounded, still keeping to two variables,
then the resulting binary partition operations does not necessarily reduce to
one of the sixteen--due to the complicated compounding of the closure
operations. For instance, the Pierce's Law formula $\left(  \left(
\sigma\Rightarrow\pi\right)  \Rightarrow\sigma\right)  \Rightarrow\sigma$ for
partitions is not equivalent to the constant function $\mathbf{1}$ since it is
not a partition tautology. The topic of the total number of binary operations
on partitions obtained by compounding the sixteen basic binary Boolean
operations is one of many topics in partition logic that awaits future research.

\section{Concluding Remarks}

In conclusion, perhaps some remarks are in order as to why it took so long to
extend the Boolean operations to partitions. The Boolean operations are
normally associated with subsets of a set or, more specifically, with
propositions. Boole originally defined his logic as the logic of subsets
\cite{boole:lot} of a universe set. It is then a theorem that the same set of
subset tautologies is obtained as the truth-table tautologies. Perhaps because
\textquotedblleft logic\textquotedblright\ has been historically associated
with propositions, the texts in mathematical logic throughout the twentieth
century (to the author's knowledge) ignored the Boolean logic of subsets and
started with the special case of the logic of propositions and then took the
truth-table characterization as the \textit{definition} of a tautology.

By the middle of the twentieth century, category theory was defined
\cite{eilenberg-mac:natequi} and the category-theoretic duality was
established between subobjects and quotient objects, e.g., between subsets of
$U$ and quotient sets (or equivalently equivalence relations or partitions) of
$U$. The conceptual cost of restricting subset logic to the special case of
propositional logic is that subsets have the category-theoretic dual concept
of partitions while propositions have no such dual concept. Hence the focus on
\textquotedblleft propositional logic\textquotedblright\ did not lead to the
search for the dual logic of partitions (\cite{ell:lop},
\cite{ell:intropartlogic}) or to the simple and natural application of Boolean
operations to partitions as well as subsets--which has been our topic here.


\begin{thebibliography}{9}                                                                                                %


\bibitem {Birk:lt}Birkhoff, Garrett 1973. \textit{Lattice Theory (}$3^{rd}%
$\textit{\ Ed.)}. New York: American Mathematical Society.

\bibitem {boole:lot}Boole, George. 1854. \textit{An Investigation of the Laws
of Thought on Which Are Founded the Mathematical Theories of Logic and
Probabilities}. Cambridge: Macmillan and Co.

\bibitem {bmp:eqrel}Britz, Thomas, Matteo Mainetti and Luigi Pezzoli 2001.
Some operations on the family of equivalence relations. In \textit{Algebraic
Combinatorics and Computer Science: A Tribute to Gian-Carlo Rota}. H. Crapo
and D. Senato eds., Milano: Springer: 445-59.

\bibitem {eilenberg-mac:natequi}Eilenberg, Samuel, and Saunders Mac Lane.
1945. General Theory of Natural Equivalences. \textit{Transactions of the
American Mathematical Society} 58 (2): 231--94.

\bibitem {ell:lop}Ellerman, David 2010. The Logic of Partitions: Introduction
to the Dual of the Logic of Subsets. \textit{Review of Symbolic Logic}. 3 (2
June): 287-350.

\bibitem {ell:intropartlogic}Ellerman, David 2014. An Introduction to
Partition Logic. \textit{Logic Journal of the IGPL.} 22, no. 1: 94--125.

\bibitem {rota:logiccers}Finberg, David, Matteo Mainetti and Gian-Carlo Rota
1996. The Logic of Commuting Equivalence Relations. In \textit{Logic and
Algebra}. Aldo Ursini and Paolo Agliano ed., New York: Marcel Dekker: 69-96.

\bibitem {ore:ter}Ore, Oystein 1942. Theory of equivalence relations.
\textit{Duke Mathematical Journal}. 9: 573-627.

\bibitem {sachs:iden}Sachs, David. 1961. Identities in Finite Partition
Lattices. \textit{Proceedings of the American Mathematical Society} 12 (6): 944--45.
\end{thebibliography}
\end{document}